\documentclass[12pt]{amsart}
\usepackage{amssymb}
\usepackage{amsmath}
\usepackage{amsbsy}
\usepackage{amscd}
\advance\textheight by \topskip
\makeatletter
\def\cal{\mathcal}
\def\Bbb{\mathbb}

\newenvironment{pf*}[1]{\proof[#1]}{\endproof}
\newcommand{\rom}{\textup}
\renewcommand{\thesubsection}{\thesection(\@roman\c@subsection)}
\makeatother
\newtheorem{Theorem}[equation]{Theorem}

\newtheorem{Lemma}[equation]{Lemma}
\newtheorem{Proposition}[equation]{Proposition}

\theoremstyle{definition}

\newtheorem{Example}[equation]{Example}

\makeatletter
\renewcommand\section{\@startsection{section}{1}%
  {\z@}{.7\linespacing\@plus\linespacing}{.5\linespacing}%
  {\reset@font\normalfont\bfseries\centering}}
\makeatother

\theoremstyle{remark}
\newtheorem{Remark}[equation]{Remark}

\newtheorem*{Acknowledgements}{Acknowledgements}
\numberwithin{equation}{section}
\numberwithin{figure}{section}

\newcommand{\thmref}[1]{Theorem~\ref{#1}}

\newcommand{\lemref}[1]{Lemma~\ref{#1}}
\newcommand{\propref}[1]{Proposition~\ref{#1}}

\newcommand{\Romnum}[1]{\expandafter\uppercase\expandafter{\romannumeral #1}} 
\newcommand{\C}{{\Bbb C}}
\newcommand{\Z}{{\Bbb Z}}

\newcommand{\R}{{\Bbb R}}

\newcommand{\CP}{\operatorname{\C P}}

\newcommand{\U}{\operatorname{\rm U}}
\newcommand{\SO}{\operatorname{\rm SO}}
\newcommand{\Spin}{\operatorname{\rm Spin}}
\newcommand{\Spinc}{\Spin^{c}}
\newcommand{\STop}{\operatorname{\rm STop}}
\newcommand{\SpinTop}{\operatorname{\rm SpinTop}}

\newcommand{\Map}{\operatorname{Map}}

\newcommand{\Aut}{\operatorname{Aut}}

\newcommand{\sign}{\operatorname{sign}}

\newcommand{\M}{{\cal M}} 
\newcommand{\G}{{\cal G}} 

\newcommand{\Met}{\operatorname{Met}}
\newcommand{\Diff}{\operatorname{Diff}}

\newcommand{\Sc}{{\cal S}}
\newcommand{\I}{{\cal N}}

\newcommand{\id}{\operatorname{id}}

\newcommand{\X}{\Bbb X}
\begin{document}
\title[Smoothability of $\Z\times\Z$-actions on 4-manifolds]{Smoothability of $\Z\times\Z$-actions on 4-manifolds}
\author{Nobuhiro Nakamura}
\address{Graduate School of Mathematical Sciences, University of Tokyo, 3-8-1, Komaba, Meguro-ku, Tokyo, 153-8914, Japan}
\email{nobuhiro@ms.u-tokyo.ac.jp}
%
\begin{abstract}
We construct a nonsmoothable $\Z\times\Z$-action on the connected sum of an Enriques surface and $S^2\times S^2$, such that each of the generators is smoothable.
We also construct a nonsmoothable self-homeomorphism on an Enriques surface.
\end{abstract}
\maketitle
%
%
\section{Introduction}\label{sec:intro}
%
%
The purpose of this paper is to prove the existence of a nonsmoothable $\Z\times\Z$-action on a $4$-manifold, such that each of the generators is smoothable: 
\begin{Theorem}\label{thm:main}
Let $X$ be the connected sum of an Enriques surface with $S^2\times S^2$.
Then, there exists a pair $(f_1, f_2)$ of self-homeomorphisms of $X$ which has the following properties{\rom :} 
\begin{enumerate}
\item $f_1$ and $f_2$ commute.
\item Each one of $f_1$ and  $f_2$ can be smoothed for some smooth structure on $X$.
However, $f_1$ and  $f_2$ can not be smoothed at the same time for any smooth structure on $X$.
\end{enumerate}
\end{Theorem}
We also construct a nonsmoothable self-homeomorphism of an Enriques surface.
\begin{Theorem}\label{thm:Enriques}
There exists a self-homeomorphism of an Enriques surface $Y$ which is nonsmoothable with respect to any smooth structure on $Y$.
\end{Theorem}
To prove these results, we modify the argument in \cite{Nf} which analyses the Seiberg-Witten moduli for families, and give more convenient constraints on diffeomorphisms of $4$-manifolds, and then, construct homeomorphisms which violate the constraints.
\begin{Acknowledgements}
The author would like to thank M.~Furuta, Y.~Kametani, K.~Kiyono and  M.~Ue for helpful discussions and comments on earlier versions of the paper.
It is also his pleasure to thank the referee for detailed comments and many valuable suggestions which enable him to make the paper corrected and improved.
\end{Acknowledgements}
%
%
\section{Constraints on diffeomorphisms}\label{sec:constraint}
%
%
In this section, we review the paper \cite{Nf}, and give some modifications of its results. 
In the paper \cite{Nf}, the author investigated the Seiberg-Witten moduli of families of $4$-manifolds, and as an application, gave some constraints on diffeomorphisms of $4$-manifolds.
Let $X$ be a closed oriented smooth $4$-manifold, and $B$ another closed manifold.
We assume a family $\X$ of $X$ over $B$ is given as a fiber bundle over $B$ whose fibers are diffeomorphic to  $X$ as oriented manifolds.
The fiber over $b\in B$ is denoted by $X_b$.
Let $T(\X/B)$ be the tangent bundle along the fiber of $\X$, and assume a metric on  $T(\X/B)$ is given.
In order to consider the Seiberg-Witten equations on the family $\X$, we need a family of $\Spinc$-structures on $\X$.
One can obtain such a family of $\Spinc$-structures if a $\Spinc$-structure on $T(\X/B)$ is given. 
For this purpose, we gave somewhat complicated sufficient conditions. (See Proposition 2.1 of \cite{Nf} and its correction \cite{Nf-c}.) 
In order to obtain a more convenient condition, we will take an alternative approach using classifying maps as described in \cite{Sz}.

Let $\Diff(X)$ be the group of orientation-preserving diffeomorphisms of $X$.
The classifying space $B\Diff(X)$ classifies families $\X\to B$ as above. 
Suppose a $\Spinc$-structure $c$ on $X$ is given.
Let us consider the group ${\Sc}(X,c)$ of pairs $(f,u)$, where $f$ is an orientation-preserving diffeomorphism and $u\colon f^*c\to c$ is an isomorphism.
The corresponding classifying space $B{\Sc}(X,c)$ classifies families $\X\to B$ with a $\Spinc$-structure $\tilde{c}$ on $T(\X/B)$ such that the restriction of $\tilde{c}$ to each fiber is isomorphic to $c$.
We have the forgetful map $\Phi\colon \Sc(X,c)\to\Diff(X)$. 
In general, $\Phi$ is not surjective. 
Let  $\I(X,c)$ be the image of $\Phi$. Then there is  an exact sequence
$$
1\to \G\to \Sc(X,c)\to \I(X,c)\to 1,
$$
where $\G=\Aut(c)\cong\Map(X,S^1)$.
Note that $B\I(X,c)$ classifies families $\X\to B$  whose structure groups are included in $\I(X,c)$.
The exact sequence leads to a fibration
$$
B\G\to B\Sc(X,c)\to B\I(X,c).
$$
Suppose it is given a family  $\X\to B$ classified by $\rho\colon B\to B\I(X,c)$ .
If $b_1(X)=0$, then $B\G$ is homotopic to $\CP^\infty\cong BS^1$.
In such a case, there is the sole obstruction to lift $\rho\colon B\to B\I(X,c)$ to $\tilde{\rho}\colon B\to B\Sc(X,c)$ in $H^3(B;\Z)$.
In particular, if $\dim B\leq 2$, then every $\rho\colon B\to B\I(X,c)$ has a lift  $\tilde{\rho}\colon B\to B\Sc(X,c)$.

Two kinds of families whose structure groups are in $\I(X,c)$ will be used in the proofs of propositions below.
The first is a mapping torus $X_f=(X\times [0,1])/f\to S^1$ defined by a diffeomorphism $f\colon X\to X$ satisfying $f^*c\cong c$.
The second is a ``double'' mapping torus $X_{(f_1,f_2)}\to S^1\times S^1$ defined by two commutative diffeomorphisms $f_1$ and $f_2$ satisfying $f_1^* c\cong f_2^*c\cong c$.
If  the family $\X$ is $X_f$ or $X_{(f_1,f_2)}$ as above, we always have a $\Spinc$-structure on $T(\X/B)$ by the previous paragraph.

When a $\Spinc$-structure $\tilde{c}$ on $T(\X/B)$ is given, the Seiberg-Witten moduli space for the family $\X$ is given as follows.
Let us define the bundle of parameters $\Pi\to B$ by
$$
\Pi = \{(g_b,\mu_b)\in\Met(X_b)\times\Omega^2(X_b)\,|\,*_b\mu_b=\mu_b\},
$$
where $\Met(X_b)$ is the space of Riemannian metrics on $X_b$ and $*_b$ is the Hodge star for the metric $g_b$.
If we choose a section $\eta$ of $\Pi$, then the moduli space for the family $(\X,\tilde{c})$ is defined by
$$
\M(\X,\tilde{c},\eta)=\coprod_{b\in B}\M(X_b,c_b,\eta_b),
$$
where $\M(X_b,c_b,\eta_b)$ is the Seiberg-Witten moduli space of the fiber $X_b$ with the $\Spinc$-structure $c_b=\tilde{c}|_{X_b}$ for the parameter $\eta_b=(g_b,\mu_b)$. 

With these understood, we can modify the results in \cite{Nf} as follows.
For a $\Spinc$-structure $c$ on $X$, let $L$ be the determinant line bundle of $c$.
Then the virtual dimension $d(c)$ of the Seiberg-Witten moduli of $(X,c)$ is given by,
$$
d(c) =\frac14(c_1(L)^2-\sign(X)) - (1+b_+).
$$
\begin{Proposition}\label{prop:prop1}
Let $X$ be a closed oriented smooth $4$-manifold with $b_1=0$ and $b_+=1$, $c$ a $\Spinc$-structure on $X$ with $d(c)=0$, and $f\colon X\to X$ an orientation-preserving diffeomorphism. 
If $f^*c$ is isomorphic to $c$, then $f$ preserves the orientation of $H^+(X;\R)$.
\end{Proposition}

The proof of \propref{prop:prop1} is given by a slight modification of the proof of Theorem 1.2 of \cite{Nf}. 
For reader's convenience, we outline it briefly.
Suppose a diffeomorphism $f$ satisfying $f^*c\cong c$ is given, and consider the mapping torus $X_f\to B=S^1$ by $f$.
Under the assumptions of \propref{prop:prop1}, the moduli space $\M(X_f,\tilde{c},\eta)$ of $X_f$ for a generic choice of $\eta$ is a compact $1$-dimensional manifold whose boundary points consist of reducibles.
Let us introduce a vector bundle $H^+_\eta\to B$ by $H^+_\eta = \coprod_{b\in B} H^+_{g_b}$, where $H^+_{g_b}$ is the space of $g_b$-self-dual harmonic $2$-forms of $X_b$. 
(Such a bundle $H_\eta$ is defined for any family.)
Then, it is proved the number of reducibles is equal modulo $2$ to the number of zeros of a generic section of $H^+_\eta$.
If $f$ reverses the orientation of $H^+(X;\R)$, then $H^+_\eta$ is a nontrivial real line bundle over $S^1$.
This is a contradiction, because the number of boundary points of a compact $1$-dimensional manifold is even.

Similarly, we can prove the following by modifying the proof of Theorem 1.1 of \cite{Nf}:
\begin{Proposition}\label{prop:prop2}
Let $X$ be a closed oriented smooth $4$-manifold with $b_1=0$ and $b_+=2$, and $c$ a $\Spinc$-structure on $X$ with $d(c)=-1$.
Suppose a pair $(f_1,f_2)$ of orientation preserving diffeomorphisms on $X$ satisfies the following conditions{\rom :}
\begin{enumerate}
\item $f_1$ and $f_2$ commute. 
\item $f_1$ and $f_2$ preserve the isomorphism class of $c$.
\end{enumerate}
Then, $w_2\left(H^+_{\eta}\right)=0$, where $H^+_{\eta}$ is the bundle associated to $X_{(f_1,f_2)}$.
\end{Proposition}
%
%
%
%
\section{Nonsmoothable self-homeomorphism on Enriques surface}\label{sec:Enriques}
%
%
The purpose of this section is to prove \thmref{thm:Enriques}.
First, note that the Enriques surface can be decomposed into three connected summands {\it topologically} by a theorem due to Hambleton and Kreck\cite{HK}.
In fact, the following theorem can be proved from Theorem 3 in \cite{HK} and its proof.
\begin{Theorem}[Hambleton-Kreck~\cite{HK}]\label{thm:HK}
The Enriques surface is homeomorphic to a topological manifold $Y=|E_8|\# \Sigma\#(S^2\times S^2)$, where $|E_8|$ is the ``$E_8$-manifold'', i.e., the simply-connected closed topological $4$-manifold whose intersection form is the negative definite $E_8$, and $\Sigma$ is a non-spin rational homology $4$-sphere with fundamental group $\Z/2$.
\end{Theorem}
\begin{Remark}
Neither $\Sigma$ nor $|E_8|\#(S^2\times S^2)$ is smoothable, because both have nontrivial Kirby-Siebenmann invariants.
\end{Remark}
Now, we will construct a self-homeomorphism of $Y$.
Let $\varphi\colon S^2\times S^2\to  S^2\times S^2$ be an orientation-preserving diffeomorphism which has the following properties:
\begin{enumerate}
\item There is a $4$-ball $B_0\subset  S^2\times S^2$ such that the restriction of  $\varphi$ to $B_0$ is the identity map on $B_0$. 
\item $\varphi$ reverses the orientation of $H^+(S^2\times S^2;\R)$.
\end{enumerate}
Such a $\varphi$ can be easily constructed as follows:
\begin{Example}\label{ex:involution}
Assume $S^2\times S^2 = \CP^1\times\CP^1$. Let $\varphi_0$ be the automorphism on $\CP^1\times\CP^1$ defined by the complex conjugation.
Choose a fixed point $p_0$ of $\varphi_0$.
Then, a required $\varphi$ is obtained by perturbing $\varphi_0$ around $p_0$ to be the identity on a neighborhood of $p_0$.
\end{Example}

Let us define a self-homeomorphism $f$ on $Y$ by $f=\id_{|E_8|\#\Sigma}\#\varphi$, where $\id_{|E_8|\#\Sigma}$ is the identity map of $|E_8|\#\Sigma$.
(Note that we can take a connected sum of $\varphi$ with $\id_{|E_8|\#\Sigma}$ on $B_0\subset S^2\times S^2$.)
Now, we claim that $f$ is nonsmoothable with respect to any smooth structure on $Y$.

To prove $f$ nonsmoothable, we will temporarily need a {\it topological} $\,\Spinc$-structure on the topological manifold $Y$.
Let us make a digression for it.
(A brief description for topological spin structures is found in \cite{Edmonds}, Section 3. See also \cite{FQ}, 10.2B.)
By Kister-Mazur's theorem, the tangent microbundle $\tau Y$ determines up to isomorphism the topological ``frame'' bundle $F$ whose structure group $\STop(4)$ consists of orientation-preserving homeomorphisms of $\R^4$ preserving the origin. 
It is known that the inclusion $\SO(4)\to \STop(4)$ induces an isomorphism on $\pi_1$ and both have trivial $\pi_0$ and $\pi_2$ (\cite{KS},V and \cite{FQ}, 8.7).
Let $\phi \colon \SpinTop(4)\to \STop(4)$ be the unique double covering.
Then, a topological spin structure on $Y$ is defined as a double covering $\tilde{F}\to F$ whose restriction to each fiber is $\phi$. 
Topological $\Spinc$-structures are similarly defined by using $\SpinTop^c(4) :=\SpinTop(4)\times_{\Z_2}\U(1)\to \STop(4)$.
The set of isomorphic classes of topological $\Spinc$-structures has a principal action of $H^2(Y;\Z)$ as in the case of true $\Spinc$-structures.
\begin{Lemma}\label{lem:pres}
Let $c$ be the topological $\Spinc$-structure on $Y$ whose $c_1(L)$ is a torsion class.
Then $f^*c$ is isomorphic to $c$.
\end{Lemma}
\proof
In this proof, all spin/$\Spinc$-structures are understood as topological ones. 
The $\Spinc$-structure $c$ can be identified with the sum of the unique spin structure $c_0$ on $|E_8|\# (S^2\times S^2)$ and a $\Spinc$-structure $c_\Sigma$ on $\Sigma$ whose  $c_1(L)$ is a torsion class.
Since $f$ is the identity on $\Sigma$, $f$ preserves $c_\Sigma$.
On the other hand, since $c_0$ is the unique spin structure on $|E_8|\# (S^2\times S^2)$, $f^*c_0\cong c_0$.
\endproof
Let us prove $f$ nonsmoothable.
Once a smooth structure on $Y$ is given, we have a reduction of the topological frame bundle $F$ to the true frame $\SO(4)$-bundle, and also a topological $\Spinc$-structure is reduced to the corresponding true $\Spinc$-structure.
Suppose $f$ is smoothed. 
By \lemref{lem:pres}, $f^*c$ is isomorphic to $c$ as true $\Spinc$-structures.
On the other hand, $f$ is an orientation-preserving diffeomorphism which reverses the orientation of $H^+(Y)$. 
This contradicts  \propref{prop:prop1}.
%
%
%
\section{Proof of \thmref{thm:main}}\label{sec:comm}
%
%
%
In this section, we prove \thmref{thm:main}.
To begin with, we collect the ingredients needed for our construction.
Let $S_0=S^2\times S^2$, and fix a $4$-ball $B_0^\prime\subset S_0$.
For $i=1,2$, let $(S_i,\varphi_i)$ be copies of $(S^2\times S^2, \varphi)$, and  fix smooth $4$-balls $B_i^\prime\subset S_i$ on which $\varphi_i|_{B_i^\prime}$ are the identity maps.
If we make a connected sum of $S_i$ ($i=0,1,2$) with another manifold, remove $B_i^\prime$ from $S_i$ and glue it along the boundary to another.
Let $Z$ be $|E_8|\# \Sigma$. 
Later, we will choose $4$-balls $B_0$, $B_1$ and $B_2$ in $Z$ so that
\begin{itemize}
\item $B_1\cap B_0=\emptyset$, $B_1\cap B_2=\emptyset$, and
\item if we make a connected sum of $Z$ with $S_i$ ($i=0,1,2$), remove $B_i$ from $Z$ and glue $\overline{Z\setminus B_i}$ and $\overline{S_i\setminus B_i^\prime}$. 
(The resulting connected sum will be denoted as $Z\#_{\partial B_i} S_i$.)
\end{itemize}
Let $E_1$ and $E_2$ be smooth $4$-manifolds homeomorphic to an Enriques surface.
The basic idea of our construction is as follows.
The connected sum $S_1\#_{\partial B_1}Z\#_{\partial B_2}S_2$ can be assumed as a connected sum of an Enriques surface with $S^2\times S^2$ in two ways: $S_1\#E_1$ and $E_2\#S_2$.
Then, commutative two homeomorphisms $f_1$, $f_2$ will be defined by $\varphi_1\#\id_{E_1}$ and $\id_{E_2}\#\varphi_2$,

Let us begin the precise construction.
Choose a $4$-ball $B_0\subset Z$ arbitrarily. 
Then $Z\#_{\partial B_0}S_0$ is homeomorphic to an Enriques surface.
Fix a homeomorphism $\bar{h}_1\colon E_1\to Z\#_{\partial B_0}S_0$.
Next, choose $B_1$ so that $D_1:=\bar{h}_1^{-1}(B_1)$ is a {\it smoothly embedded $4$-ball} in $E_1$. 
Take a smooth connected sum $S_1\#_{\partial D_1}E_1$ and a (topological) connected sum $S_1\#_{\partial B_1}Z\#_{\partial B_0}S_0$ so that a homeomorphism $h_1=\id_{S_1}\#\bar{h}_1\colon S_1\#_{\partial D_1}E_1\to S_1\#_{\partial B_1}Z\#_{\partial B_0}S_0$ is defined. 

Note that $S_1\#_{\partial B_1}Z$ is also homeomorphic to an Enriques surface.
Fix a homeomorphism $\bar{h}_2\colon E_2\to S_1\#_{\partial B_1}Z$.
Choose $B_2$ so that $D_2:=\bar{h}_2^{-1}(B_2)$ is a {\it smoothly embedded $4$-ball} in $E_2$. 
Take a smooth connected sum $E_2\#_{\partial D_2}S_2$ and a (topological) connected sum $S_1\#_{\partial B_1}Z\#_{\partial B_2}S_2$ so that a homeomorphism $h_2=\bar{h}_2\#\id_{S_2}\colon E_2\#_{\partial D_2}S_2\to S_1\#_{\partial B_1}Z\#_{\partial B_2}S_2$ is defined.

Define the self-diffeomorphism $\bar{f}_1$ on $S_1\#_{\partial D_1}E_1$ by $\bar{f}_1=\varphi_1\#\id_{E_1}$, and $\bar{f}_2$ on $E_2\#_{\partial D_2} S_2$ by $\bar{f}_2=\id_{E_2}\#\varphi_2$.
Choose a homeomorphism  $h\colon S_1\#_{\partial B_1}Z\#_{\partial B_2}S_2 \to S_1\#_{\partial B_1}Z\#_{\partial B_0}S_0$ so that $h|_{S_1\setminus B_1^\prime}$ is the identity map.
Via homeomorphisms $h$, $h_1$ and $h_2$, we obtain self-homeomorphisms $f_1$ and $f_2$ of $X:=S_1\#_{\partial B_1}Z\#_{\partial B_2}S_2$ induced from $\bar{f}_1$ and $\bar{f}_2$, respectively.
Then each $f_i$ ($i=1,2$) is smoothable for the smooth structure $E_i\#_{\partial D_i} S_i$.
Clearly, $f_1$ and $f_2$ commute. 
Let $c$ be the $\Spinc$-structure on $X$ whose $c_1(L)$ is a torsion class.
As in \lemref{lem:pres}, we can see that $f_1$ and $f_2$ preserve the isomorphism class of $c$.
However, $w_2\left(H^+_{\eta}\right)\neq 0$  by construction.
By \propref{prop:prop2}, $f_1$ and $f_2$ can not be smoothed at the same time.
Thus, \thmref{thm:main} is proved.
%
%
%
\section{Remarks}\label{sec:remarks}
%
%
%
We give two remarks.
The first is on another possibility of application of \propref{prop:prop2}.
The following problem would be interesting: {\it Find two diffeomorphisms of a smooth manifold homeomorphic to a connected sum of an Enriques surface $E$ with $S^2\times S^2$ that are simultaneously smoothable, commute up to isotopy, but do not have representatives in their isotopy classes that commute. }
If we want to construct such two diffeomorphisms on the {\it smooth} connected sum $E\#S^2\times S^2$, then one of the difficulties would be as follows. 
To appeal to \propref{prop:prop2}, one of two diffeomorphisms will be required to reverse the orientation of the $H^+(E)$-part of $H^2(E\#S^2\times S^2)$, and it will be easy if we can construct such a diffeomorphism as a connected sum of a diffeomorphism $f$ of $E$ with one of $S^2\times S^2$. 
However, this method is impossible, because \propref{prop:prop1} prohibits such an $f$.

The second remark is on a generalization of the construction of the moduli spaces for families. 
In fact, we can construct the moduli space for a family without a family of $\Spinc$-structures.
More precisely, we claim the following: {\it When a family $\X\to B$ is classified by  $\rho \colon B\to B\I(X,c)$,  we can always construct the moduli space $\M(\X,c)$ for the family $\X$, even if $\rho$ does not have a lift  $\tilde{\rho}\colon B\to B\Sc(X,c)$.}
The construction is outlined as follows.
By taking local trivializations, the family $\X$ can be given via transition functions $\psi_{\beta\alpha}\colon U_\alpha\cap U_\beta \to \I(X,c)$ for an appropriate covering $\{U_\lambda\}_{\lambda\in\Lambda}$ of $X$.
Suppose the intersection of every two members in $\{U_\lambda\}_{\lambda\in\Lambda}$ is contractible.
Then we can take a lift of each $\psi_{\beta\alpha}\colon U_\alpha\cap U_\beta \to \I(X,c)$ to $\tilde{\psi}_{\beta\alpha}\colon U_\alpha\cap U_\beta \to \Sc(X,c)$. In general, such $\tilde{\psi}_{\beta\alpha}$ do not satisfy the cocycle condition, but satisfy it {\it up to gauge}, i.e., $\psi_{\gamma\beta}\psi_{\beta\alpha}\psi_{\gamma\alpha}^{-1}$ is a gauge transformation.
One can define local families $\M(U_\lambda\times X,c)=\coprod_{b\in U_\lambda}\M(\{b\}\times X,c)\to U_\lambda$ of moduli spaces and attaching maps $\tilde{\psi}_{\beta\alpha}^*$ between them induced from  $\tilde{\psi}_{\beta\alpha}$. 
(Here, we need a little care on metrics and perturbations.)
Since the moduli spaces are defined as the quotient spaces divided by the gauge transformations,  $\tilde{\psi}_{\beta\alpha}^*$ satisfy the cocycle condition.
Therefore, the global family $\M(\X,c)$ can be constructed from the local families $\M(U_\lambda\times X,c)$ via $\tilde{\psi}_{\beta\alpha}^*$.
Such a family $\M(\X,c)$ would be useful for further applications.


\begin{thebibliography}{99}


\bibitem{Edmonds} A.~Edmonds, 
{\it Aspects of group actions on four-manifolds},
Topology Appl. {\bf 31} (1989), no. 2, 109--124. 

\bibitem{FQ} M.~H.~Freedman and F.~Quinn, 
Topology of 4-manifolds,
Princeton Mathematical Series, 39. Princeton University Press, Princeton, NJ, 1990. 


\bibitem{HK} I.~Hambleton and M.~Kreck, 
{\it Smooth structures on algebraic surfaces with cyclic fundamental group},
Invent. Math. {\bf 91} (1988), no. 1, 53--59. 

\bibitem{KS} R.~C.~Kirby and L.~C.~Siebenmann,
Foundational essays on topological manifolds, smoothings, and triangulations.
Annals of Mathematics Studies, No. 88. Princeton University Press, Princeton, N.J.; University of Tokyo Press, Tokyo, 1977.


\bibitem{Nf} N.~Nakamura, 
{\it The Seiberg-Witten equations for families and diffeomorphisms of 4-manifolds},
Asian J. Math. {\bf 7} (2003), no. 1, 133--138; 

\bibitem{Nf-c} \bysame,
Correction, Asian J. Math. {\bf 9}  (2005),  no. 2, 185. 

%
%

\bibitem{Sz} M.~Szymik,
{\it Characteristic cohomotopy classes for families of 4-manifolds},
preprint.

\end{thebibliography}
\end{document}